\documentclass[12pt]{amsart}

\usepackage{amssymb,graphicx}
\usepackage[all]{xy}
 \usepackage[T1]{fontenc}

\numberwithin{equation}{section}

\newcommand{\be}{\begin{enumerate}}
\newcommand{\ee}{\end{enumerate}}
\newcommand{\bi}{\begin{itemize}}
\newcommand{\ei}{\end{itemize}}

\def\R{\mathbb{R}}

\def\al{\alpha}
\def\be{\beta}

\def\De{\Delta}

\def\si{\sigma}

\def\ep{\varepsilon}

\def\nd{\noindent}
\def\bull{\hfill$\Box$\\}
\def\proof{\nd {\bf Proof.\ }}

\def\p{\partial}

\begin{document}

\today
\vskip 1cm
\begin{center}
{\sc  About the Bloch-Connelly-Henderson Theorem on the simplexwise linear homeomorphisms 
of a convex 2-disk

\vspace{1cm}

Jean Cerf}

\end{center}

\title{}
\author{}
\address{15, rue Sarrette, 75014 Paris, France}
\email{Jean Cerf <jean.cerf@math.u-psud.fr>}


\begin{abstract}
 This paper gives an improved version of the original proof 
of  Bloch-Connelly-Henderson's theorem about the space of $SL$ homeomorphisms of a convex 2-disk.
A major improvement is related to the main lemma of the original paper.
\end{abstract}
\maketitle

\thispagestyle{empty}
\vskip 1cm
\section*{Introduction}

A simplexwise linear disk, briefly an $SL$ \emph{disk}, of $\R^2$  is a submanifold with boundary homeomorphic to the 2-disk $D^2$ and triangulated by means 
of affine simplices. Its boundary is called an $SL$ \emph{circle} of $\R^ 2$. An $SL$ circle of $\R^2$ 
is said to be \emph{convex} if the disk it bounds is convex, and \emph{strictly convex} if in addition it is never flat at a 
vertex. The 1984 paper by Bloch, Connelly and Henderson (that will be quoted as [BCH]) proves the following results.\\

\nd{\bf Theorem 1.} {\it  Let $K$ be an $SL$ disk of $\R^2$ and $f$ an $SL$ embedding of $\p K$ in $\R^2$. If $f(\p K)$ 
is strictly convex, 
or more generally if $K$ has no $f$-obstructive 1-simplex {\rm (see \S 4)} then $f$ can be extended to an $SL$ 
embedding of $K$ in $\R^2$.}\\

\nd{\bf Theorem 2.} {\it If $K$ is convex and has $m$ interior vertices the space $\mathcal E_{\p K}(K)$ of $SL$ homeomorphisms of $K$ keeping $\p K$ fixed, endowed with the compact-open topology, is homeomorphic to $\R^{2m}$.} \\

Theorem 2 was stated as the main result in the introduction of [BCH]. Although Theorem 1 was not stated, it is implicitely contained in the final statement ([BCH], Theorem 5).
This paper aims to give a simplified and clarified version of the original proofs. It is self-contained, giving 
detailed proofs even of the parts whose treatment is similar to the original one (\S 1 and \S2).

What makes the proof of these theorems difficult is the fact that $\mathcal E_{\p K}(K)$ is not an invariant of the abstract 
model of $K$. In other words, two $SL$ disks $K$ and $K'$ may be $SL$ isomorphic while the corresponding spaces 
$\mathcal E_{\p K}(K)$ and $\mathcal E_{\p K'}(K')$ are not homeomorphic.

Clearly, $\mathcal E_{\p K}(K)$ is a $2m$-dimensional manifold; indeed, ordering the set of the interior vertices of 
$K$ identifies $\mathcal E_{\p K}(K)$ with an open subspace of $\R^{2m}$. The study of that manifold easily
reduces to the case where $K$ is strictly convex and \emph{simple} ({\it i.e.} no 1-simplex separates 
$K$ in two parts).  The proof proceeds by induction on the number $p$ of 2-simplices of $K$.

A projective transformation of the plane  
changes neither convexity nor simplicity of $K$ nor again the topology of $\mathcal E_{\p K}(K)$. So, one may suppose that $K$
has been put in  \emph{reduced form}. This means that some 1-simplex $\si$ of $\p K$ is the $[0,1]$ segment of the $x$-axis, the remaining arc of  $\p K$ being the graph of a concave $SL$ function over this segment (see Lemma 1).
 The 1-simplex $\si$
 is a facet of some 2-simplex $\tau$ whose third vertex $s$ is interior. Removing
$\tau$ from $K$ leads therefore to an $SL$ disk $L$ with $(p-1)$ 2-simplices. The disk $L$ is not convex, but deforming
$\p L $ into a strictly convex $SL$ circle is easy; indeed, keeping all vertices except $s$ fixed, one moves $s$  vertically({\it i.e.} along a parallel to the $y$-axis)
 to a point $s'$ outside of $K$ and close to $\si$. In order to extend  this deformation 
 to the whole $L$, the authors
made use of what they called the ``basic lemma'' 
([BCH], Lemma 3.1).

By introducing \emph{keys} and \emph{twin-keys}, Lemma 3 of the present paper leads to a much simpler version
of the ``basic lemma'' (Lemma 4), which is stated in the strictly convex case as:

{\it If the $SL$ disk $K$ is transverse to the verticals ({\it i.e.} any intersection with a parallel of the $y$-axis is connected), then any vertical $SL$ embedding $f:\p K\to \R^2$ with strictly 
convex image can be extended to a vertical $SL$ embedding of $K$ ({\it i.e.} an embedding 
moving vertically each vertex). 
} 

The  extension from $\p L$ to $L$ 
completes
the inductive step in the proof of Theorem 1. 
Concerning Theorem 2, it remains to prove that $\mathcal E_{\p L}(L)$ is invariant under this extension. 
The present version of this part of the proof takes advantage of the fact that Lemma 4 works directly in the space of  $SL$ embeddings whereas the ``basic lemma'' works in the space of  $SL$ mappings with non-negative volumes.  But the general line essentially remains the original one.
Let us say that $f$ and $f'$ in $\mathcal E_{\p L}(L)$
are equivalent if their composition with the projection of $\R^2$ to the $x$-axis coincide (in other words, if $f'\circ f^{-1}$
is vertical). The proof relies on the fact that this equivalence relation leads to a product decomposition 
with fibers homeomorphic to $\R^m$.\\

I cannot close this introduction wihout expressing my gratitude to Fran\c cois Laudenbach. 
Without his constant and amical interest and help---including both mathematical and practical points of view---this paper would not exist.\\

\nd{\sc Conventions and notations.}
The following conventions hold for the whole paper. The coordinates in $\R^2$ are denoted $(x,y)$. The projection 
to the $x$-axis (resp. $y$-axis) is denoted by $\pi_x$ (resp. $\pi_y$). The lines parallels to
$Ox$ (resp. $Oy$) are called horizontals (resp. verticals).

The coordinates in $\R^n\times \R^{n'}$ are denoted by $\left((x_1,\ldots,x_n),(y_1,\ldots,y_{n'})\right)$
or by $(X,Y)$. The two projections are respectively denoted $\Pi_x$ and $\Pi_y$. 

The vertices of an $SL$ disk being denoted by $s_1,\ldots,s_n$, the space of all its $SL$ mappings to $\R^2$
is identified with $\R^n\times \R^n$ by the homeomorphism
$$f\mapsto \left( (\pi_x f(s_1),\ldots, \pi_xf(s_n)),(\pi_y f(s_1),\ldots, \pi_yf(s_n))\right).
$$
\medskip
\section{Projective reduction of a convex $SL$ circle of $\R^2$}

Let $C$ be a convex $SL$ circle of $\R^2$. A maximal line segment 
in $C$ will be called a \emph{natural edge}.  
The $SL$ circle 
whose 1-simplices are the natural edges of $C$ 
 is called the \emph{strictly convex $SL$ circle
associated with $C$}. 

A projective transformation of $\R^2$ maps $C$ to an $SL$ circle if and only if the inverse image of the infinity line 
does not meet $C$; 
convexity and strict convexity are preserved.\\
\begin{figure}[h]
  ${}$\quad \quad\includegraphics[scale=.4]{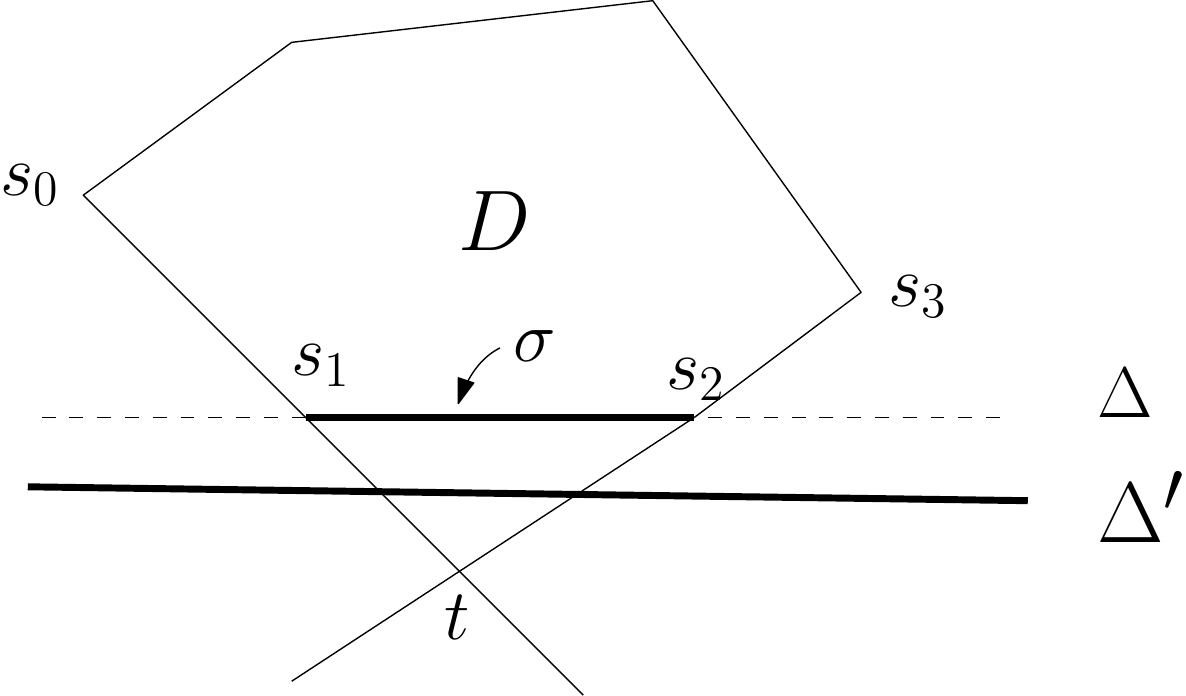}
 \hfill
  \includegraphics[scale=.4]{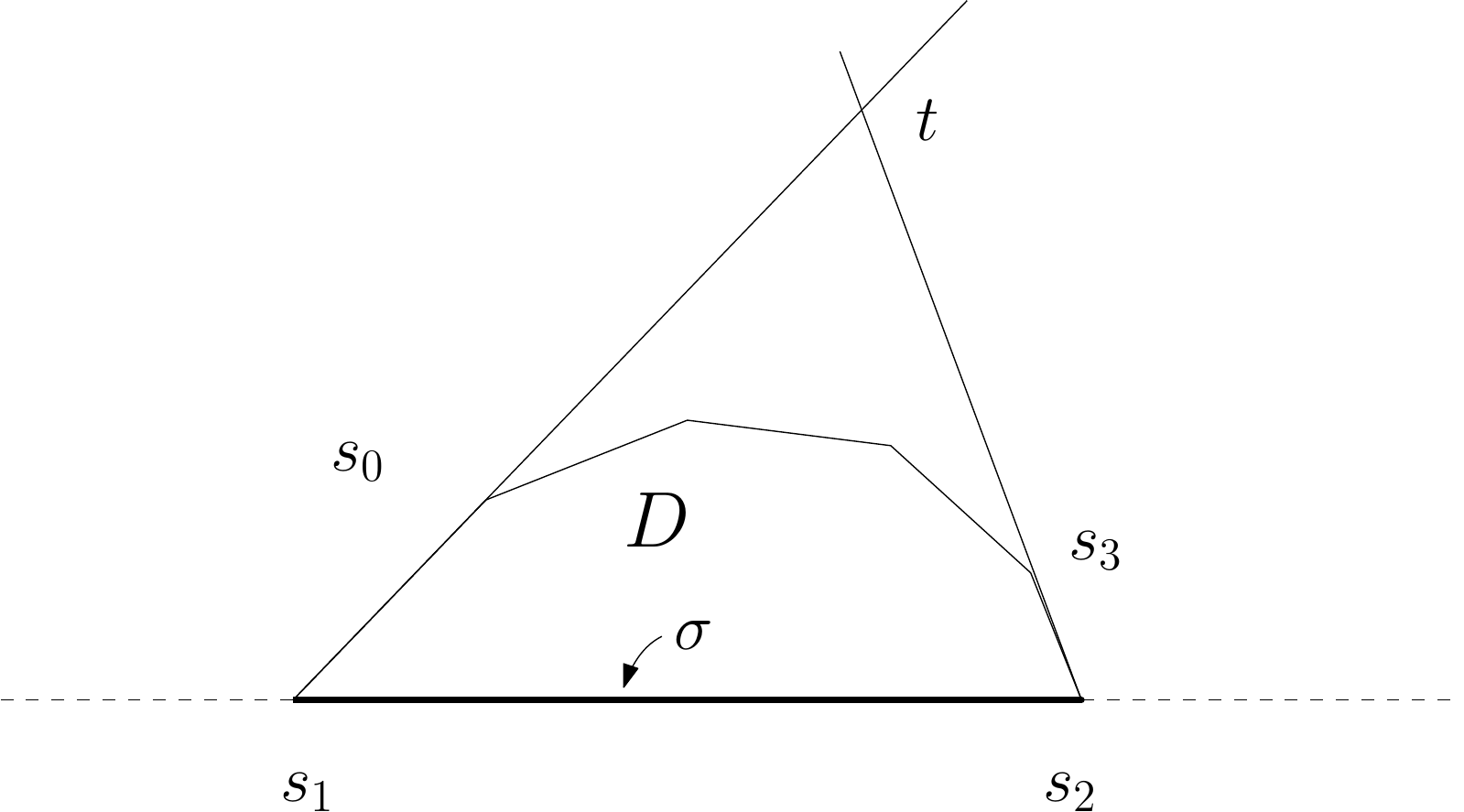}\quad\quad
 \caption{}\label{cerf2}
 \end{figure}

\nd{\bf Lemma 1.} {\it Let $C$ be a convex $SL$ circle and $\mu$ be a natural edge of $C$. There exists a projective 
transformation $g$ of $\R^2$ such that $g(C)$ is an $SL$ circle, $g(\mu)$ is $[0,1]\times \{0\}$ and 
$g(C\smallsetminus\mu)  \subset (0,1)\times \R^+$.}\\

One says that such a $g$ puts $C$ in \emph{reduced form}. When $C$ is in reduced form and bounds  an 
$SL$ disk $K$, then $K$ is also said to be in reduced form. \\

\proof The general case easily reduces to the strictly convex one: one replaces $C$ by the associated strictly convex
$SL$ circle. 

Let $\si$ be any 1-simplex of the strictly convex $C$. Let $s_1$ and $s_2$ be the vertices of $\si$ and 
let $s_0$ (resp. $s_3$) be the vertex right before $s_1$ (resp. right after $s_2$). One denotes by $D$ the disk
bounded by $C$.
The line containing $[s_0,s_1]$ and the line containing $[s_2, s_3] $  are distinct since $D$
is salient at $s_1$ and $s_2$. So, they meet at a unique point $t$ (maybe at infinity) which is not on the line 
$\De$ carrying $\si$. The case of Figure \ref{cerf2} left ($t$ and $D$ are on opposite sides of $\De$)
and the case where $t$ is at infinity, both cases  reduce to the case of Figure \ref{cerf2} right (both $t$ and $D$
 are at finite distance and lie on the same side of $\De$) by means of a projective transformation of $\R^2$
 whose inverse image $\De'$ of  the line at infinity intersects the open 1-simplices $(t,s_1)$ and $(t, s_2)$.
 In that case, $u$  being any interior point of $\si$, the disk $D$, except $s_1$ and $s_2$,
 is contained in the open  band 
 limited by the two parallels to the direction $\overrightarrow{ut}$ passing respectively through $s_1$ and $s_2$.  
 The proof is achieved 
 by means of a suitable affine isomorphism of $\R^2$. \bull 
 \section{$TrV$ disks: embeddings and mappings with non-negative volumes}
 
 Let $K$ be an $SL$ disk of $\R^2 $ whose vertices are denoted by $s_i \ (i=1,\ldots, n)$
 and the 2-simplices by $\tau_j \ (j=1, \ldots, p)$. The space of all mappings from $K$ to 
 $\R^2$ being identified to $\R^n\times \R^n$, the function which maps 
 $f$ to the algebraic volume of $f(\tau_j)$, noted $vol(f(\tau_j))$ (see Section 7), identifies to
  an alternate bilinear form $\ell_j$
 on the vector space $\R^n$
 taking a positive value at $f$ when $f:K\hookrightarrow\R^2$ is the natural injection.
 The space of all $SL$ mappings with non-negative volumes identifies with the convex
 polyhedral cone defined by 
 $$ (1)\quad\quad \ell_j(X,Y)\geq 0 \quad \text{for }j= 1,\ldots, p\,. $$
 
 The system obtained by adding to (1) the equation $X=X_0$, with $X_0=\Pi_x(s_1,\ldots, s_n)$
 defines the subspace of vertical $SL$ mappings. We will usually write 
 $$(2)\quad\quad \ell_{j,X_0}(Y)\geq 0 \quad\text{for } j=1,\ldots p\,,
 $$
 where $\ell_{j,X_0}$ is a linear form on $\{X_0\}\times \R^n$ identified to $\R^n$. This linear form has a positive value at 
 the point $\Pi_y(s_1,\ldots, s_n)$.\\
 
 \nd{\bf Definitions 1.} Let $K$ be an $SL$ disk of $\R^2$. One says that $K$
 is \emph{transverse to the verticals} (we will say $TrV)$) if any non-empty
 intersection of $K$ with a vertical line is connected. 
 
In this setting, the \emph{roof} of $K$ is the 1-disk of $\p K$ made of the 1-simplices
$\si$ such that no ascending vertical half-line issued from a point of $\si$ meets the interior of $K$.

Similarly, one defines a $TrH$ disk (transverse to the horizontals).\\

\nd{\bf Remark.} Some 1-simplices of the roof of a $TrV$ disk, including the first and the last ones, may be vertical.\\

\nd {\bf Notations.} Let $K$ be a $TrV$ disk. One denotes by $\mathcal W(K)$ the space of vertical 
$SL$ mappings from $K$ to $\R^2$ with non-negative volumes (recall that vertical means that each vertex moves 
vertically). One denotes  by $\mathcal W_T(K)$ the subspace 
of mappings which fix the roof $T$ of $K$. One denotes by $\mathcal V(K)$ (resp. $\mathcal V_T(K)$) the subspace
of $\mathcal W(K)$ (resp. $\mathcal W_T(K)$) consisting of the orientation preserving embeddings.

The vertices of $K$ being numbered from 1 to $n$ one assumes that the first $q$ ones belong to $T$.\\

\nd{\bf Lemma 2.} {\it The space $\mathcal W(K)$  identifies to 
the convex polyhedral cone of dimension $n$ of $\R^n$
defined by the system $(2)$. The subspace $\mathcal W_T(K)$  identifies to
 the convex polyhedron of dimension
$n-q$ of $\R^{n-q}$ defined by $(2)$ and 
$$(3)\quad\quad (x_i,y_i)= (\pi_x(s_i),\pi_y(s_i))\quad\text{for } i=1, \ldots, q\,.$$
The space $\mathcal V(K)$ is the interior of $\mathcal W(K)$ and $\mathcal V_T(K)$ is the interior 
of $\mathcal W_T(K)$.}\\

\proof Clearly $\mathcal V$ (resp. $\mathcal V_T$) is a non-empty open subset of $\mathcal W$ (resp. $\mathcal W_T$).
It remains to prove that $\mathcal V$ is the interior of $\mathcal W$, in other words that any 
$f$ which is solution of the system
$$(2')\quad\quad \ell_{j,X_0}(Y)>0 \quad\text{for } j=1,\ldots , p$$
is an injective mapping. A sufficient condition for this is that $f\vert (K\cap V)$ is injective on every vertical line $V$.
Now, $K\cap V$ is either one point or, $K$ being $TrV$,  is an $SL$ 1-disk. By $(2')$ $f$ is strictly increasing
on any 1-simplex of this 1-disk.\bull

\section{$TrV$ disks continued: keys and twin-keys} 

\nd{\bf Definitions 2.} Let $K$ be a $TrV$ disk. A 2-simplex $\tau$ of $K$ is said to be a \emph{key}
if one of its faces $\si$ is in the roof $T$
and if the \emph{link} of  $\si$ ({i.e.} the vertex of $\tau$ opposite to $\si$) is interior to $K$ 
and projects vertically to the interior of $\si$.

A pair $(\tau_\ell, \tau_r)$ of adjacent 2-simplices is said to be a \emph{twin-key} if $\tau_\ell\cap T$ and $\tau_r\cap T$
are 1-simplices whose common link is interior to $K$ and if in addition $\tau_\ell\cap\tau_r$ is vertical.\\

\nd{\bf Definitions 3.} Let $K$ be an $SL$ disk of $\R^2$.

An  1-simplex $\si$ of $K$ is said to be \emph{spanning} if  
$\p \si=\si \cap  \p K$. 

  If $K$ has no spanning 1-simplex, 
it is said to be \emph{simple}.\\

\nd{\bf Remark.} Definitions 3 depend only on the abstract model of $K$.\\

\nd {\bf Lemma 3.} {\it Let $K$ be a simple $TrV$ disk having more than one 2-simplex. Then $K$ possesses a key or a twin-key.}\\

\nd{\bf Remark.} A simple $SL$ disk having more than one 2-simplex has at least three 2-simplices.\\

\proof  One denotes by $s_i \ (i=1, \ldots,q)$ from the left to the right the vertices of the roof, by $\si_i\ (i= 1,\ldots, q-1)$
the 1-simplex starting at $s_i$, by $\tau_i$ the 2-simplex whose $\si_i$ is a facet, and by $s'_i$ the link of $\si_i$.
Since $K$ is simple, $s'_i\in \p K$ would imply $\p \tau_i\subset\p K$ and hence $K= \tau_i$ which is excluded.
So, $s'_i$ is an interior point of $K$ for every $i$. 

If $q=2$, then $\tau_1$ is a key. So, we may suppose $q>2$.
If $\tau_1$ is not a key then $\pi_x(s'_1)\geq \pi_x(s_2)$. Suppose  that some $i\leq q-2$ satisfies 
$\pi_x(s'_i)\geq \pi_x(s_{i+1}$); if neither $\si_{i+1}$ is a key nor $(\si_i,\si_{i+1}) $ is a twin-key, then 
$\pi_x(s'_{i+1})\geq \pi_x(s_{i+2})$. 
This is impossible if $\si_{i+1}$ is vertical descending or if $i=q-2$.\bull
\vskip -2cm
 \section{$TrV$ disks completed: the main Lemma}
 
 \nd{\bf Definition 4.} Let $K$ be an $SL$ disk
 and $f$ be an $SL$ embedding of $\p K$ in $\R^2$ with a convex image. A spanning 1-simplex $\si$
 is said to be $f$-\emph{obstructive} if the $f$-image of one of the arcs that $\p \si$ bounds in $\p K$ is flat.\\
 
 \nd{\bf Property 1.} {\it Let $\si$ be a spanning 1-simplex which divides $K$ into $K_1$ and $K_2$. 
 If $\si$ is not $f$-obstructive then $f$ canonically defines an $SL$ embedding $f_i: \p K_i\to \R^2$.
 It has the property  that 
 if $K$ has no $f$-obstructive spanning 1-simplex, then $K_i$ is so for $f_i$ $(i= 1,2)$.}\\
 
 \proof If $K_i$ has a spanning 1-simplex $\si_i$, then the arc of $\p K_i$ containing $\si$ cannot be flattened 
 by $f$ because $f(K)$ is convex. The same is true for the remaining arc because $\si_i$ is not $f$-obstructive.\bull\\
 
 \nd{\bf Property 2.} One supposes that $K$ is $TrV$ with roof $T$ and that $f$ is vertical. If $K$
 has an $f$-obstructive spanning 1-simplex $\si$ with one end point in $T\smallsetminus \p T$ then the other 
 end point lies in $T$.\\
 
 \proof Suppose $\si$ has one end point in $\p K\smallsetminus T$. Then both arcs of $\p K\smallsetminus \p \si$
 would have non-injective projection on the $x$-axis. So, they could not be flattened by any vertical embedding. 
 
 ${}$\bull
 
  \begin{figure}[h]
 ${}$ \quad\quad\includegraphics[scale =.5]{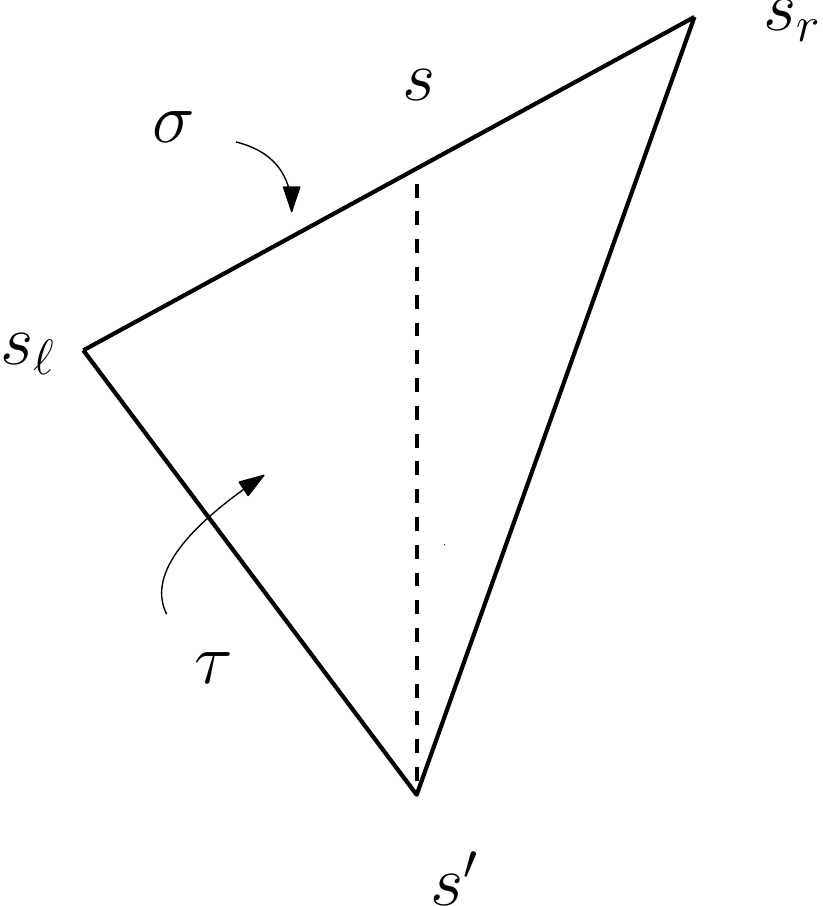}
\hskip 2cm
  \includegraphics[scale =.5]{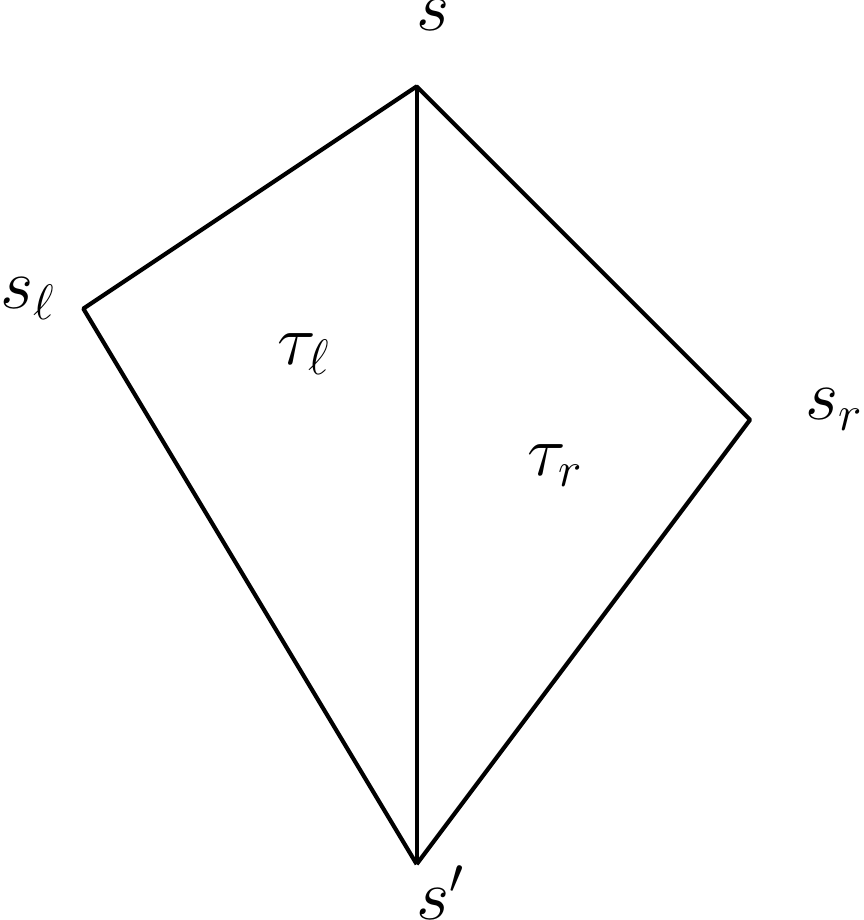} \quad\quad
 \caption{On the left, a \emph{key};  on the  right, a \emph{twin-key}.}\label{cerf4}
 \end{figure}
 
 The notations concerning keys and twin-keys that we will use in the proof of the next lemma are respectively shown on Figures \ref{cerf4} left and right. On Figure  \ref{cerf4} left, $s$ is the vertical projection of $s'$ to $\si$.\\

\nd{\bf Lemma 4 (the main lemma).} {\it Let $K$ be a $TrV$ disk and $v$ be an $SL$ vertical embedding
of $\p K$ to $\R^2$. If $v(\p K)$ is convex and $K$ has no $v$-obstructive 1-simplex then $v$ extends 
to an $SL$ vertical embedding of $K$ to $\R^2$.\\
}

\proof One may suppose that $v$ has the positive orientation. The method is an induction on the number $p$
of 2-simplices of $K$. The statement is obvious for $p=1$. So, one assumes that $p>1$ and that the statement is proved
for every $TrV$ disk having a number of 2-simplices less than $p$.

\begin{figure}[h]
  \includegraphics[scale =.5]{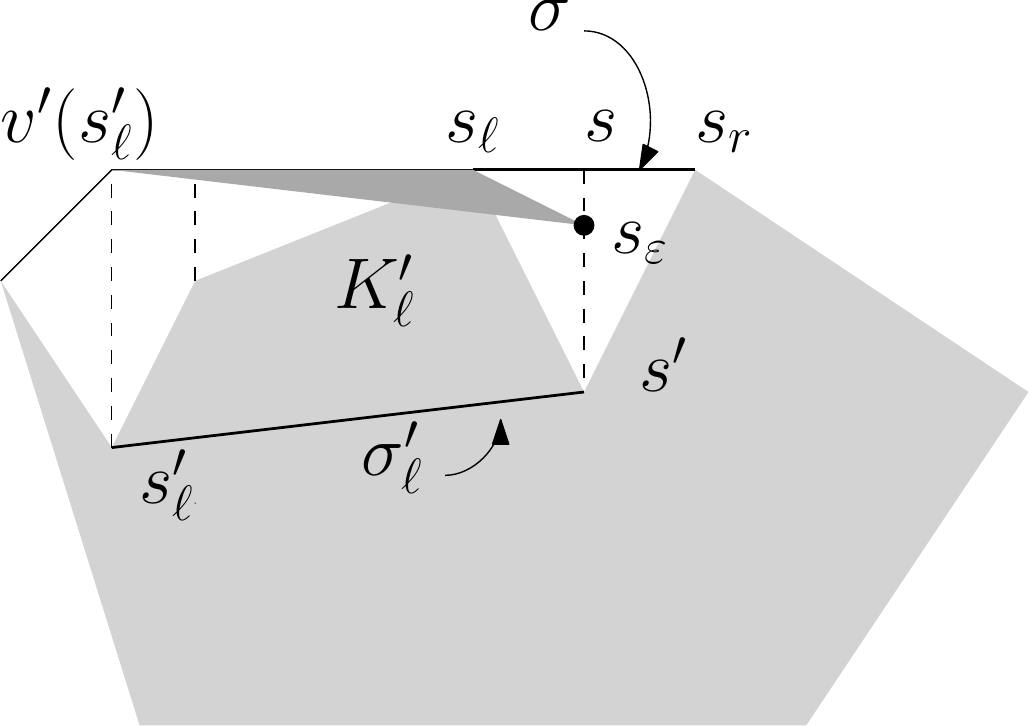}
 \caption{Here is the special case where $K'_r$ is empty and a part of $\p K$
 including $\si$ is kept fixed by $v$. The dark grey region is $v_\ep(K'_\ell)$; the light grey region is  $K'$.}\label{cerf5}
 \end{figure}

If $K$ has a spanning 1-simplex $\si$ which divides $K$ into $K_1$ and $K_2$, then by Property 1 both $K_1$
and $K_2$ satisfy the induction assumption. The $SL$ vertical embeddings defined by $v$ on $\p K_1$ 
and $\p K_2$ coincide along $\si$. Therefore, their vertical extensions to $K_1$ and $K_2$ define a vertical extension
of $v$ to $K$.
So, we may suppose that $K$ is 
simple
, and hence, has a key $\tau$ or a twin-key 
$\tau_\ell\cup\tau_r $ (se Lemma 3 and Figure 
\ref{cerf4}). Recall that in both cases the vertex $s'$ 
is interior to $K$. One denotes by $K'$ (with roof $T'$) the $TrV$ disk obtained from $K$ by removing $\tau$
(this case is represented on Figure  \ref{cerf5}) 
or by removing $\tau_\ell\cup\tau_r$ in case of a twin-key.

Let $v'$ be the $SL$ vertical embedding $\p K'\to \R^2$ which coincides with $v$ on $\p K'\cap\p K$ and 
verifies $v'(s')=v(s)$. As $K$ is simple $s'$ belongs to $\p \si'$ for every spanning 1-simplex $\si'$ of $K'$.
If $\si'$ is $v'$-obsructive, Property 2 implies that the other vertex $s''$ of $\si'$ lies in $T'$; so, $\si'$ cannot be vertical.  
  Assume for instance that $s''$ lies to the left of $s'$ then $\si'$ is $v'$-obstructive if and only if $v'(s'')$ is aligned 
  with $v(\si)$ in case of a key and with $v(\si_\ell)$ in case of a twin-key.
  Denote by $s'_\ell$ the leftmost point $s''$ when $\si'$ runs over all $v'$-obstructive spanning 1-simplices
  of $K'$; let $\si'_\ell$ be the corresponding spanning 1-simplex. 
  Let $K'_\ell$ be the part of $K'$ over
  $\si'_\ell$. There are similar definitions of $\si'_r$ and $K'_r$.
  One denotes by $K''$ the $TrV$ disk obtained from $K'$ by removing $K'_\ell\cup K'_r$; let $v'':\p K''\to \R^2$
  be the $SL$ vertical embedding which coincides with $v'$ on $\p K''\cap \p K'$. The induction assumption
  applies to $(K'', v'')$. So, $v''$ extends to an $SL$ vertical embedding of $K''$. 
  Such an extension still exists if one replaces $v''$ by the $SL$ vertical embedding $v''_\ep$ 
 which coincides with $v''$ at every vertex of  $\p K''$ except at $s'$ whose image moves to a point $s_\ep$ close to $s$ on its descending vertical half-line.
   The disk $K'_\ell$ endowed 
  with the $SL$ vertical embedding $v'_\ell$ of $\p K'_\ell$ which coincides with $v$ on $\p K'_\ell\cap \p K$
  and verifies $v'_\ell(s')= v(s_\ep)$ satisfies the induction assumption. The same is  true for $K'_r$ and also for $\tau$
  (or $\tau_\ell\cup \tau_r$).
  
  The corresponding extensions, as well as the extension of $v''_\ep$ to $K''$, all given by the induction assumption,
  coincide two by two on the intersections of their respective domains and coincide with $v$ 
  on the intersections of these domains with $\p K$. Gluing altogether these extensions gives an extension of $v$
  to $K$.\bull
 \vskip -.8cm
 
 \nd{\bf Corollary.} {\it Let $K$ be a $TrV$ disk. Assume that the roof $T$ is concave and $K$ has no 1-simplex
 $\si$ joining the end points of $T$;
   suppose also that $\p K\smallsetminus T$ has only one vertex $u$. Then the image 
 of the evaluation map at $u$, $v\mapsto v(u)$, 
 of the space of all $SL$ vertical embeddings of $K$ into $\R^2 $ fixing $T$ is a vertical open half-line, infinite in the descending direction, whose upper bound lies above the line segment joining the end points of $T$.\\
}

\proof We may suppose that the end points of $T$ are the points 0 and 1 of the $x$-axis. Denote the projection
$\pi_x(u)$ by $u'$. One has $0<u'<1$. If $\pi_y(v(u))<0$, the image $v(\p K)$ is strictly convex and hence,  there
exists no $v$-obstructive spanning 1-simplex. This property remains true if $v(u)=u'$, because a 
$v$-obstructive spanning 1-simplex should have one end point at $u'$, which is impossible. The fact that the evaluation mapping has an open image completes the proof.

\section{Proof of Theorem 1}

Recall the statement of Theorem 1: $K$ is an $SL$ disk of $\R^2$, $f:\p K\to\R^2$ is an $SL$ embedding.
If $K$ has no $f$-obstructive spanning 1-simplex then $f$ extends to an $SL$ embedding $K\to \R^2$.\\

\nd{\sc Special case: $f(\p K)$ is strictly convex.}

The proof proceeds by induction of the number $p$ of 2-simplices of $K$. The case $p=1$ is obvious. One assumes
that the result is proved up to $p-1$ with $p> 1$. If $K$ has a spanning 1-simplex dividing $K$ into $K_1$ and 
$K_2$, both are strictly convex and $f$ defines an embedding $f_i: \p K_i\to \R^2 \ (i= 1,2)$. By induction assumption,
$f_i$ extends to an $SL$ embedding $K_i\to\R^2$. Gluing these extensions together yields an extension
of $f$ to $K$. So, one may suppose that $K$ is simple. Moreover by Lemma 1, one may assume that $f(\p K)$
is in reduced form. 

\begin{figure}[h]
${}$\quad\quad\includegraphics[scale = .6]{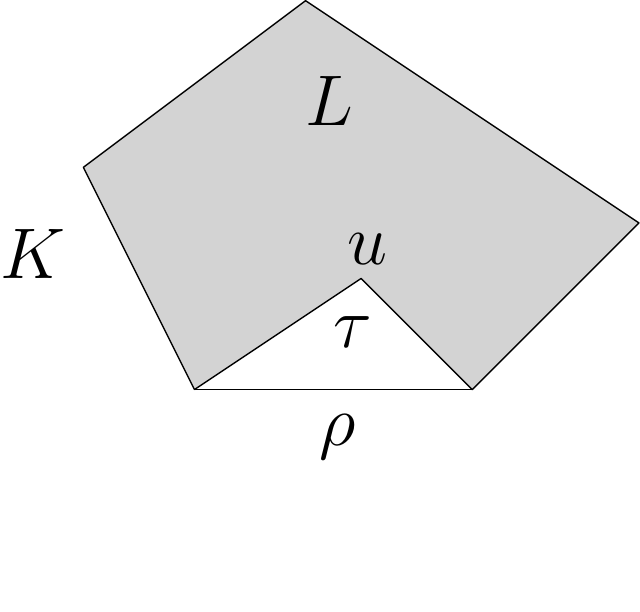}
 \quad \quad\quad
  \includegraphics[scale =.5]{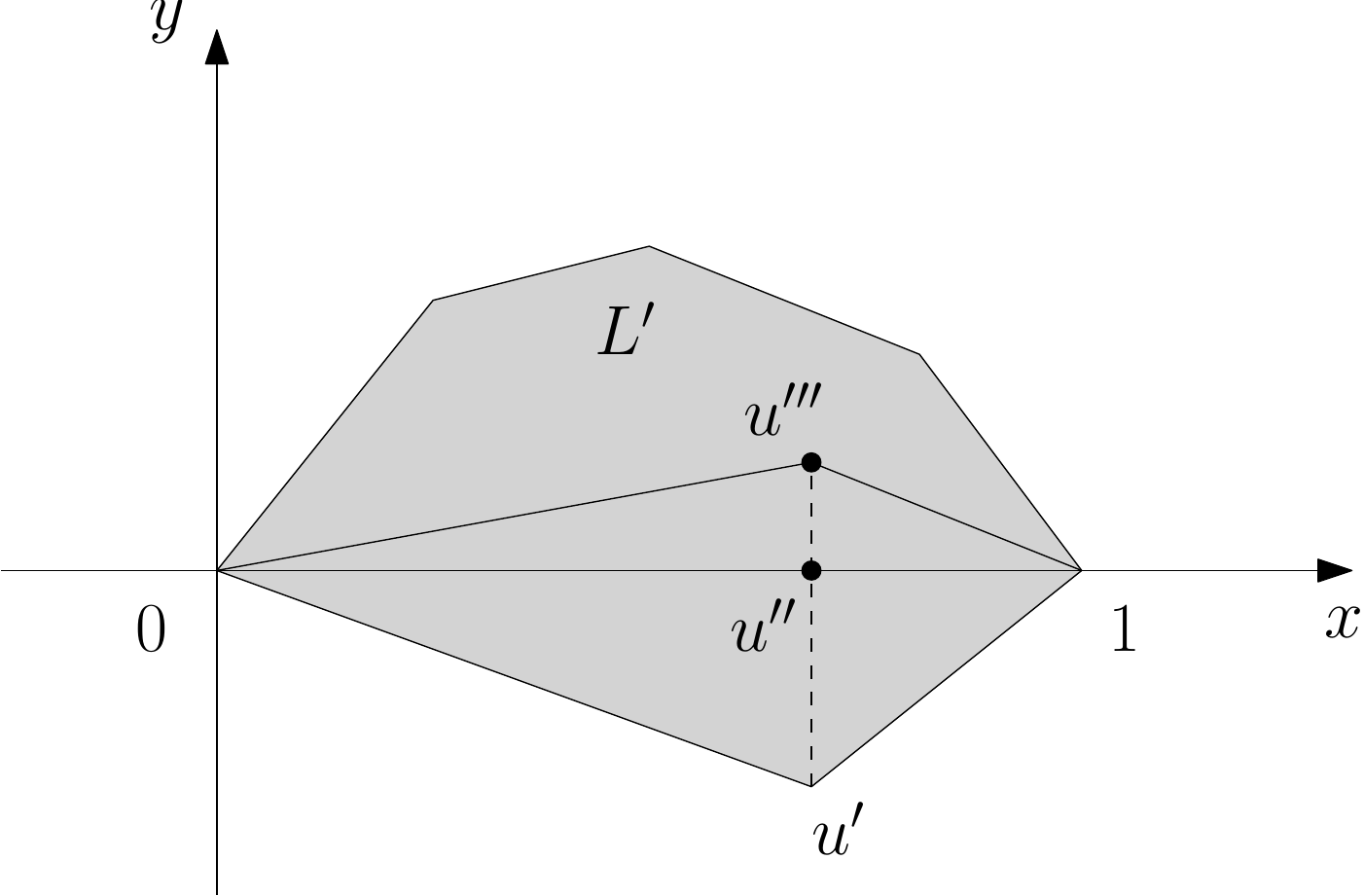}\quad
 \caption{}
 \label{cerf7}
 \end{figure}

Let $\rho$ be the 1-simplex of $\p K$ such that $f(\rho)=[0,1]\times\{0\}$;  one denotes by $\tau$ the 2-simplex of 
$K$ whose $\rho$ is a face and by $u$ the link of $\rho$. As $K$ is simple, 
$u$ is an interior vertex of $K$. One denotes by $L$ the $SL$ disk obtained by removing $\tau$ from $K$ (Figure \ref{cerf7}).
 \medskip

Let $u'$ be a point of $\R^2$ below $(0,1)\times\{0\}$. Denote by $f'$ the $SL$ embedding $\p L\to\R^2$
which coincides with $f$ on $\p L\cap \p K$ and maps $u$ to $u'$. The image $f'(L)$ is strictly convex and, by
induction, there exists an $SL$ embedding extending $f'$ to $L$. One still denotes this extension by $f'$ and one denotes its image
by $L'$ (see Figure \ref{cerf7}); set $u''=\pi_x(u')$ and choose a point $u'''$ on the vertical ascending from $u''$. Let
$v$ be the $SL$ embedding $\p L'\to\R^2$ fixing the roof and sending $u'$ to $u''$. 
Applying Lemma 4 to $(L',v)$ 
yields an extension $f''$ of $v$ to $L'$. Denote by $f'''$
the $SL$ mapping $L'\to\R^2$ which maps $u'$ to $u'''$ and fixes the other vertices. If $u'''$ is close enough to
$u''$ then $f'''$ is an embedding. By construction, $f'''\circ f'$ has a canonical extension to $K$.\bull

\nd{\sc General case.} 

\begin{figure}[h]
  \includegraphics[scale =.5]{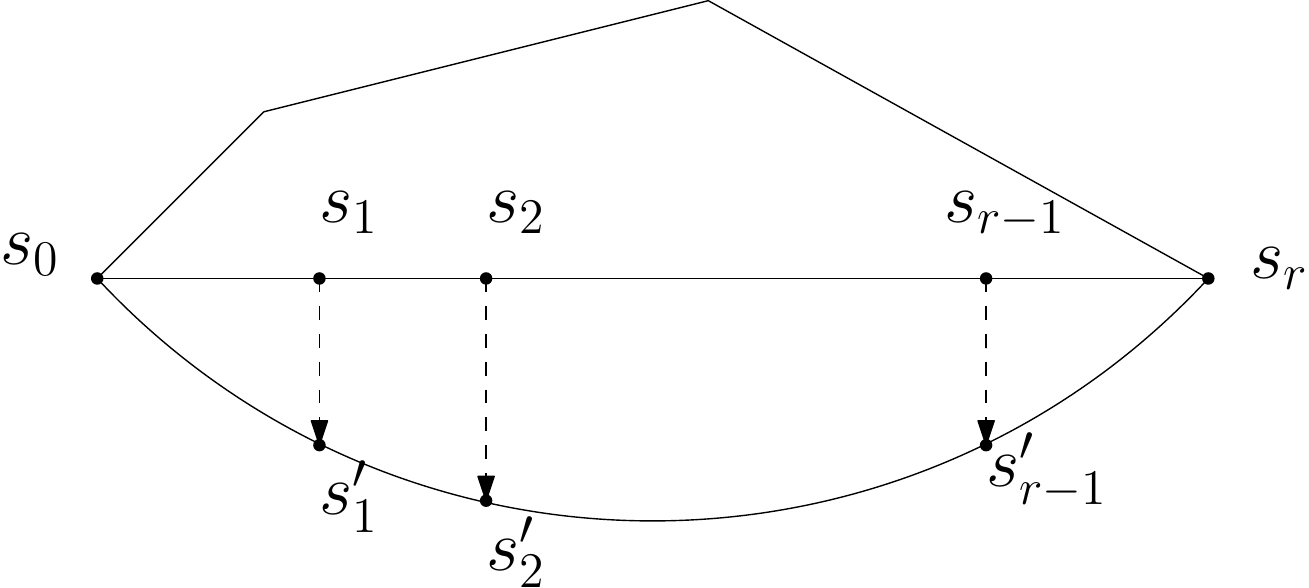}
 \caption{}\label{cerf8}
 \end{figure}

Let us name a \emph{plateau} a natural edge with more than one 1-simplex.
The previous case may be called the ``zero plateau'' case. From it, one argues by induction on the number of plateaus.
One still supposes $f(\p K)$ in reduced form. One denotes by $s_0, s_1,\ldots, s_r$ the vertices of $f(\p K)$
on the $x$-axis. Let $A$
 be the graph of a strictly convex function on $[0,1]$ with end points $s_0$ and $s_r$, for instance the one 
 yielded by a suitable circle with center $(1/2, y), \ y>0$. Let $s'_i, \ i= 1,\ldots, r-1$, be the intersection of $A$
 with the descending vertical of $s_i$ (see Figure \ref{cerf8}). Denote by $v$ the $SL$ embedding $f(\p K)\hookrightarrow\R^2$ 
 sending $s_i$ to $s'_i$ for $i=1,\ldots, r-1$, and fixing the other vertices. Replacing $f$ by 
 $v\circ f$ removes one plateau and does not introduce any obstructive spanning 1-simplex. Lemma 4
implies that $v^{-1}$ extends to an $SL$ embedding $f'':f'(K)\to\R^2$. Then $f''\circ f'$ is the desired extension.\bull
 \medskip
 
 \section{Cairns parametrization of continuous families of convex polyhedral disks}

Let $\left( \ell_{j,X}(Y)\right)_{j= 1,\ldots,p}$ be a finite family of affine forms in $\R^n$, where each of them
depends continuously of the parameter $X\in\R^m$. One denotes by $P_X$  the part of $\{X\}\times\R^n$
defined by 
$$(4)\quad\quad \ell_{j,X}(Y)\geq 0 \quad \text{for } j=1,\ldots, p.
$$
The subspace $P_X$ is a convex closed polyhedron of $\R^n$. When $P_X$ is bounded and has a non-empty interior
then $P_X$ is homeomorphic to a disk $D^n$ (this fact is a part of the above proof of Theorem 1). 
The following lemma is a parametrized version of this result. \\

\nd{\bf Lemma 5.} {\it Let $U$ be an open set in $\R^m$ such that $P_X$ is bounded and has  a non-empty
interior for every $X\in U$. Then the projection to $\R^m$ makes $\cup_{X\in U}P_X$ a trivial fiber bundle 
homeomorphic to $U\times D^n$. \\
}

\proof The projection $\R^n\times\R^m\to \R^m$ has over $U$
infinitely many sections $X\mapsto c_X$ where $c_X$ belongs to the interior of $P_X$ for every $X\in U$.
For instance, one chooses $c_X$ equal to the barycenter of $P_X$. One parametrizes $P_X$ radially,
that is, the center $O$ of $D^n$ goes to $c_X$ and for every $\al\in S^{n-1}$ the radius $[O,\al]$ 
is mapped affinely to the parallel ``radius'' $[c_X, a_{X,\al}]$ of $P_X$. The point $a_{X,\al}$
is common to some hyperplanes 
$$(4')\quad\quad \ell_{j, X}(Y)=0
$$
where $j$ belongs to some subset made of $p'<p$ elements of $\{1,\ldots, p\}$, say $\{1, \ldots, p'\}$. For
$(X',\al')$ close enough to $(X,\al)$ and $j\leq p'$, let $a_{X',\al',j}$ denote the point where the half-line starting 
at $c_{X'}$ in the direction of $\al'$ intersects the hyperplane $\{\ell_{j, X'}(Y)=0\}$. The length of the 
``radius'' $[c_{X'},a_{X',\al'}]$
is  the smallest distance between $c_{X'}$ and $a_{X',\al',j}$ when $j$ runs over $\{1, \ldots,p'\}$, and hence, it varies continuously.
\bull

\medskip

\section{Proof of Theorem 2}

Recall the statement of Theorem 2: $K$ is a convex $SL$ disk with $m$ interior vertices. Then the space 
$\mathcal E_{\p K}(K)$ of $SL$ homeomorphisms of $K$ keeping $\p K$ fixed is homeomorphic to $\R^{2m}$.

The main part of the proof deals with the strictly convex case; the general case will easily follow.
As for Theorem 1, the method is by induction on the number of 2-simplices of $K$. 
Assume the result is proved if the number of $2$-simplices is less that $p$. One may suppose that $K$
is in reduced form and simple (the proof of this claim is left to the reader).

One keeps the notations $\tau, \ L, \ T$ of \S 5. One denotes by $s_1,\ldots, s_m$ the interior vertices of $K$
so that $s_m$ is the vertex of $\tau$ interior to $K$. One will use the following notations:
\begin{itemize}
\item $\mathcal E_T(L)$ is the space of $SL$ embeddings with positive orientation, $L\hookrightarrow \R^2$,
fixing $T$ and such that $\pi_xf(s_m)\in (0,1)$;
\item $\mathcal F_T(L)$ is the space of $SL$ mappings with non-negative volume, $L\to\R^2$, fixing $T$
and such that $\pi_xf(s_m)\in [0,1]$;
\item $\mathcal E_T^y(L), \ \mathcal E_T^{>y}(L), \ \mathcal F_T^y(L)$ etc. for which the exponent means the
part of $\R$ assigned to $\pi_yf(s_m)$.
\end{itemize}
These spaces are all identified to subspaces of $\R^n\times\R^n$. \\

\nd{\bf Lemma 6.} ${}$
\nd 1) {\it For every $y\leq 0$ the following holds:
$$ (5)\quad\quad \Pi_x\mathcal E_T(L)= \Pi_x\mathcal E^y_T(L)= \Pi_x\mathcal E^{\geq y}_T(L)= 
\Pi_x\mathcal E^{>y}_T(L).
$$
}

\nd 2) {\it For every $X\in \Pi_x\mathcal E_T(L)$ the intersections of $\Pi_x^{-1}(X)$ 
with $\mathcal F_T(L),\ \mathcal F^{\geq y}_T(L),\ \mathcal F^y_T(L)$ are convex polyhedra with respective dimension 
$m,\,  m,\,  m-1$. They are bounded except the first one. Their respective interiors are 
$\mathcal E_T(L),\ \mathcal E^{>y}_T(L),\ \mathcal E^y_T(L)$.\\
}

\proof The first item easily follows from the corollary of Lemma 4. 

For the second item,  the result concerning the pair $\left(\mathcal F_T(L), \mathcal E_T(L)\right)$
implies the other ones by intersection with half-planes or hyperplanes. Let $X\in \Pi_x\mathcal E_T(L)$ 
and $f\in \Pi_x^{-1}(X)\cap\mathcal E_T(L)$. The mapping $f'\mapsto f'\circ f^{-1}$
is an homeomorphism of the pair $\left(  \Pi_x^{-1}(X),  \Pi_x^{-1}(X)\cap \mathcal E_T(L)\right)$ onto
the pair $\left(\mathcal W_T(f(L)),\mathcal V_T(f(L)\right)$ (see \S 2 for notations). Applying Lemma 4 to the disk 
$f(L)$ completes the proof of Lemma 6. \bull 


We come back to the proof of Theorem 2. For every $y<0$ and every $f\in \mathcal E_T^y(L)$ the disk $f(L)$
is strictly convex and its number of 2-simplices is $p-1$. On the other hand $\mathcal E_{\p K}(K)$ is obviously 
homeomorphic to $\mathcal E_T^{>0}(L)$. Theorem 2 (in the stricly convex case) is therefore consequence 
of the following two homeomorphisms:
$$
\begin{array}{l}
(6)\quad\quad \mathcal E_T^{>0}(L)\cong \mathcal E_T^{-\frac 12}(L)\times (0,1),\\

(7)\quad\quad \mathcal E_T^{-\frac 12}(L)\cong \mathcal E_{\p f(L)}(f(L))\times(0,1) \quad\text{for every } f\in \mathcal E_T^{-\frac 12}(L)\,.
\end{array}
$$

\nd{\sc Proof of (6).}
Consider the space  of all $SL$ mappings $L\to \R^2$ identified with $\R^n\times\R^n$. (Notice that 
$L$ has the same vertices as $K$).
 Let $\tau_j$ be any 2-simplex of $L$ equipped 
with a numbering of its vertices which orients $\tau_j$ positively.
It defines the alternate bilinear form $vol(f(\tau_j))$
 by composing  the projection 
 $\R^n\times \R^n\to \R^3\times \R^3$ with the alternate  bilinear form on $\R^3$
 $$\left\vert\begin{array}{lll}
 1&x_1&y_1\\
 1&x_2&y_2\\
 1&x_3&y_3
 \end{array}\right\vert\ .$$
The space $\mathcal F_T(L)$ (for which only $m$ vertices are movable)
is defined in $\R^m\times\R^m$ by a system 
of $p-1$ inequations of the type (4) with the next inequation added:
$$(8) \quad\quad 0\leq x_m\leq 1\ .
$$
In order to define $\mathcal F_T^{\geq 0}(L)$ one adds to the system (4) (8)
$$(9) \quad\quad y_m\geq 0.
$$
By 2) of Lemma 6, for every $X\in \Pi_x^{-1}\left(\mathcal E_T^{>0}(L)\right)$ the set 
$\Pi_x^{-1}(X)\cap \mathcal F_T^{\geq 0}(L)$ is a convex polyhedron whose interior 
$\Pi_x^{-1}(X)\cap \mathcal E_T^{> 0}(L)$ is non-empty and bounded. 
This implies that $\Pi_x^{-1}(X)\cap \mathcal F_T^{\geq 0}(L)$ is bounded.
Therefore, Lemma 5 yields an homeomorphism
$$(10) \quad\quad \mathcal E_T^{>0}(L)\cong \left( \Pi_x\mathcal E_T^{>0}(L)\right)\times\R^m\,.
$$

The space $\mathcal F_T^{-\frac 12 }(L)$ is defined by the system (4) (8) (9)
modified by replacing $Y$ with $(y_1,\ldots, y_{m-1}, -\frac 12)$ and $y_m>0$ with 
$y_m=-\frac 12$. Applying 2) of Lemma 6 and Lemma 5 (with $n$ replaced by $m$ and $U$ replaced 
by $\Pi_x\mathcal E_T^{-\frac 12}(L)$), one gets 
$$(11)\quad\quad \mathcal E_T^{-\frac 12}(L)\cong \left(\Pi_x\mathcal E_T^{-\frac 12}(L)\right)\times\R^{m-1}\,.
$$
Applying 1) of Lemma 6 for $y=0$ and $y=-1$ yields
$$(12) \quad\quad \Pi_x\mathcal E_T^{>0}(L)\cong \Pi_x\mathcal E_T^{-\frac 12}(L)\,.
$$
Clearly (10), (11), (12) together imply (6).\bull\\

\nd{\sc Proof of (7).} Let $f\in\mathcal E_T^{-\frac 12}(L)$ with for instance $f(s_m)=(\frac 12, -\frac 12)$.
One has an homeomorphism 
$$(13)\quad\quad \mathcal E_T^{-\frac 12}(L)\cong \mathcal E_T^{-\frac 12}f((L))\,.
$$
By the chosen representation of the $SL$ mappings (see the convention at the end of the introduction),
the system which defines $\mathcal E_T^{-\frac 12}(L)$ works as well for $\mathcal E_T^{-\frac 12}f((L))$.
On the other hand, $f(L)$ is strictly convex. Therefore, it is transverse to the horizontals.
Applying the horizontal version of Lemma 4 yields
$$(14) \quad\quad \Pi_y\mathcal E_T^{-\frac 12} (f(L))\cong \Pi_y \mathcal E_{\p f(L)}(f(L))\,.
$$
Henceforth, one will see the system defining $\mathcal E_T^{-\frac 12} (f(L))$ as a system with $(m-1)$ 
parameters, namely $(y_1, \ldots, y_{m-1})$ on the $\R^m$ space of the $X$ variable.
Applying Lemma 5 in this setting, that is, $ \Pi_y\mathcal E_T^{-\frac 12} (f(L))$ as set of parameters, yields
$$(15) \quad\quad\mathcal E_T^{-\frac 12} (f(L))\cong \left(  \Pi_y\mathcal E_T^{-\frac 12} (f(L))\right)\times \R^m\,.
$$
By adding the equation $x_m=\frac 12$ in order to define $\mathcal E_{\p f(L)}(f(L))$ one  similarly obtains
$$(16) \quad\quad\mathcal E_{\p f(L)}(f(L))\cong \left(\Pi_y\mathcal E_{\p f(L)}(f(L))\right)\times\R^{m-1}\,.
$$
Clearly (13), (14), (15), (16) together imply (7).\bull\\

\nd{\sc Extension of Theorem 2 to the non-strictly convex case.}
As for theorem 1, one argues by induction on the number of plateaus. One assumes $K$ is in reduced form. 
One 
 again adopts  the notations $s_i, s'_i, v$ (see Figure \ref{cerf8}). As the extension $f'$ of $v$ to $K$ given 
by Lemma 4 is vertical the trivial fibration defined by the projection $\Pi_x$ on $\mathcal E_{\p K}(K)$
and on $\mathcal E_{\p f'(K)}(f'(K))$ have the same basis. On the other hand, both fibers are homeomorphic to 
$\R^m$. This completes the proof.

\end{document}